# Two person zero-sum fuzzy matrix games and extensions

***Rajani Singh***

*Institute of Applied Mathematics and Mechanics,*
*Universit of Warsaw ,*
*Poland*

*June 6, 2015*

---

**Abstract**

In this work, fuzzy Linear programming problems (FLPP) and their implementations is being done in two person zero-sum fuzzy matrix game theory, based on Bector, Chandra and Vijay [2] model. Verification of solution of all possible type of two person zero-sum game is done by solving examples by numerical method.

Generalized fuzzy set first introduced by Atanassov (1986) which is called I-fuzzy set has been studied and using this knowledge, application in (FLPP) and two person zero-sum matrix game with I-fuzzy goals is given. At the end I gave useful MATLAB codes for solving all three type of two person zero-sum fuzzy matrix games, with fuzzy goals, with fuzzy pay-offs and with both fuzzy goals and fuzzy pay-offs.

---

## Introduction

Similar to fuzzy linear programming problems fuzziness in matrix games can also come. In the game players may have fuzzy goals and/or the entries of the pay-off matrix are fuzzy number.

In this paper I discuss about different types of two person zero-sum fuzzy matrix game with all type of possibilities for fuzziness and their conversion in linear programing problems (LPP’s).

In real decision making problems there are instances where not only the degree to which an object belongs to a set is known but in addition a degree to which the same object does not belong to the set is also known.

Atanassov (1986) proposed a generalized approach of fuzzy set is called I-fuzzy set. I used this special type of set to solve two person zero-sum matrix game with I-fuzzy goals.

### Two person zero-sum fuzzy matrix games

There are three types of two person zero-sum matrix games listed as,

1. Two person zero-sum matrix games with fuzzy goals.
2. Two person zero-sum matrix games with fuzzy pay-offs.
3. Two person zero-sum matrix games with fuzzy goals and fuzzy Pay-offs.

### 1. Two person zero-sum matrix games with fuzzy goals:

Let $A \in R^{mxn}$ be a pay-off matrix and $S_m = \{x \in R^m : \sum_{i=1}^{m} x_i = 1, x_i \geq 0, \forall\ i=1, 2, \ldots, m\}$ and $S_n = \{x \in R^n : \sum_{j=1}^{n} x_j = 1, x_j \geq 0, \forall\ j=1, 2, \ldots, n\}$ be the strategies spaces for player I and player II respectively.

Let $v_0$ and $w_0$ be aspiration levels for player I and player II respectively and $p_0$ and $q_0$ be the tolerances for player I and player II respectively, then the two person zero-sum matrix game with fuzzy goals is denoted by FG and is defined as ,

$$FG = (S^m, S^n, A, v_0, \gtrsim, p_0, w_0, \lesssim, q_0).$$

Where $\gtrsim$ ($\lesssim$) is essentially greater (less) than or equal to sign.

**Solution of the fuzzy matrix game (FG):**

An ordered pair $(\bar{x}, \bar{y})$ is said to be a solution of the two person zero-sum matrix game with fuzzy goals FG if

$$(\bar{x})^T Ay \gtrsim_{p_0} v_0 \qquad \text{for all } y \in S^n,$$

and

$$x^T A\bar{y} \lesssim_{q_0} w_0 \qquad \text{for all } x \in S^m$$

*Conversion of fuzzy matrix game FG into fuzzy linear programming*

(FP-1) Find $x \in R^m$

such that $x^T Ay \gtrsim_{p_0} v_0$, for all $y \in S^n$

$\sum_{i=1}^{m} x_i = 1$

$x \geq 0$

(FP-2) Find $y \in R^n$

such that $x^T Ay \lesssim_{q_0} w_0$, for all $x \in S^m$

$\sum_{j=1}^{n} y_j = 1$

$y \geq 0$

Let $A_j$ (j=1,2,…n) and $A_i$ (i=1,2,…m) denote the $j^{th}$ column and $i^{th}$ row of the matrix A respectively . Then

(FP-3) Find $x \in R^m$

such that

$A_j{}^T x \gtrsim_{p_0} v_0$ , $\forall$ j=1, 2, . . . , n

$\sum_{i=1}^{m} x_i = 1$

$x \geq 0$

(FP-4) Find $y \in R^n$

such that

$A_i\, y \lesssim_{q_0} w_0$, $\forall$ i=1, 2, . . . , m

$\sum_{j=1}^{n} y_j = 1$

$y \geq 0$

Membership function for $j^{th}$ constraint of (FP-3) is given as

$$\mu_j(A_j^T x) = \begin{cases} 1 & A_j^T x \geq v_0 \\ 1 - \frac{v_0 - A_j^T x}{p_0} & (v_0 - p_0) \leq A_j^T x < v_0 \\ 0 & A_j^T x < (v_0 - p_0) \end{cases}$$

and the membership function for $i^{th}$ constraint of (FP-4) is given as,

$$\vartheta_i(A_i\, y) = \begin{cases} 1 & A_i\, y \leq w_0 \\ 1 - \frac{A_i\, y - w_0}{q_0} & w_0 < A_i\, y \leq w_0 + q_0 \\ 0 & A_i\, y > w_0 + q_0 \end{cases}$$

Then the fuzzy linear programming (FP-3) and (FP-4) become

(FLP) max $\lambda$

Subject to,

$$A_j^T x \geq v_0 - (1-\lambda)p_0 \qquad \forall\, j = 1, 2, \ldots, n$$
$$\sum_{i=1}^{m} x_i = 1$$
$$\lambda \leq 1$$
$$x, \lambda \geq 0$$

(FLD) max $\eta$

Subject to,

$$A_i y \leq w_0 + (1-\eta)q_0 \qquad \forall\, i = 1, 2, \ldots, m$$
$$\sum_{j=1}^{n} y_j = 1$$
$$\eta \leq 1$$
$$y, \eta \geq 0$$

**Remark 1.1** The crisp linear programming problems (FLP) and (FLD) are primal-dual pair in the fuzzy sense, not in crisp sense.

**Remark 1.2** If both players have the same aspiration levels i.e. $v_0 = w_0$ and in the optimal solutions of (FLP) and (FLD) $\lambda^* = \eta^* = 1$, then the fuzzy game FG reduces to crisp two person zero-sum game G.

**Example 1.1** Solve the two person zero-sum matrix game with fuzzy pay-offs FG whose pay-off matrix A is given by

$$A = \begin{bmatrix} 1 & 3 & 0 \\ 4 & 7 & 2 \\ 3 & 5 & 6 \end{bmatrix}, \quad \text{with } v_0 = 5/3,\ w_0 = 3/2 \text{ and } p_0 = 2,\ q_0 = 3.$$

**Solution.**

(FLP) max $\lambda$

Subject to;

$$2\lambda - x_1 - 4x_2 - 3x_3 \leq 1/3$$
$$2\lambda - 3x_1 - 7x_2 - 5x_3 \leq 1/3$$
$$2\lambda - 2x_2 - 6x_3 \leq 1/3$$
$$x_1 + x_2 + x_3 = 1$$
$$\lambda \leq 1$$
$$x_1, x_2, x_3, \lambda \geq 0$$

(FLD) max $\eta$

Subject to;

$$3\eta + y_1 + 3y_2 \leq 9/2$$
$$3\eta + 4y_1 + 7y_2 + 2y_3 \leq 9/2$$
$$3\eta + 3y_1 + 5y_2 + 6y_3 \leq 9/2$$
$$y_1 + y_2 + y_3 = 1$$
$$\eta \leq 1$$
$$y_1, y_2, y_3, \eta \geq 0$$

The optimal solution and optimal value for (FLP) (i.e. for player I) is $x^* = (0.0271, 0.4233, 0.5496)$ and $\lambda^* = 1.0000$ and the optimal solution and optimal value for (FLD) (i.e. for player II)

is $y^* = (0.8000, 0.0000, 0.2000)$ and $\eta^* = 0.3000$.

## 2. Two person zero-sum matrix games with Fuzzy Pay-offs:

Two person zero-sum matrix game with fuzzy pay-offs is the triplet $\widetilde{G} = (S^m, S^m, \widetilde{A})$.

**Reasonable values ($\tilde{v}, \widetilde{w}$) of the game $\widetilde{G}$:**

$(\tilde{v}, \widetilde{w})$ is called a reasonable solution of the fuzzy matrix game FG if there exists $x^* \in S^m$, $y^* \in S^n$ satisfying $(x^*)^T \widetilde{A}\, y \gtrsim \tilde{v}$, $\forall\ y \in S^n$ and $x^T \widetilde{A} y^* \lesssim \widetilde{w}$, $\forall\ x \in S^m$.

**Solution of the Fuzzy Matrix Game $\widetilde{G}$:**

T1= $\{\tilde{v} \in N(R)$: $\tilde{v}$ is reasonable value for Player I$\}$, T2=$\{\widetilde{w} \in N(R)$: $\widetilde{w}$ is reasonable value for Player II$\}$.Let there exist $\tilde{v}^* \in$ T1, $\widetilde{w}^* \in$ T2 such that $F(\tilde{v}^*) \geq F(\tilde{v})$, $\forall\ \tilde{v} \in$ T1 and $F(\widetilde{w}^*) \leq F(\widetilde{w})$, $\forall\ \widetilde{w} \in$ T2 then $(x^*, y^*, \tilde{v}^*, \widetilde{w}^*)$ is called the solution of the game $\widetilde{G}$.

**Conversion of fuzzy matrix game $\widetilde{G}$ into LPP:**

(FP) max $F(\tilde{v})$

subject to,

$$\sum_{i=1}^{m} \tilde{a}_{ij}\, x_i \gtrsim \tilde{v} - (1-\lambda)\tilde{p} \qquad \forall\ j = 1,2,\ldots,n$$
$$e^T x = 1,$$
$$x \geq 0,$$
$$0 \leq \lambda \leq 1,$$

(FD) min $F(\widetilde{w})$

subject to,

$$\sum_{j=1}^{n} \tilde{a}_{ij}\, y_j \lesssim \widetilde{w} + (1-\eta)\, \tilde{q} \qquad \forall\ i = 1,2,\ldots,m$$
$$e^T y = 1,$$
$$y \geq 0,$$
$$0 \leq \eta \leq 1,$$

By using the defuzzification function $F: N(R) \to R$ for the constraints (FP) and (FD) problems can be rewritten as:

(FP1) max $F(\tilde{v})$

subject to ,

$$\sum_{i=1}^{m} F(\tilde{a}_{ij})\, x_i \gtrsim F(\tilde{v}) - (1-\lambda)F(\tilde{p}) \qquad \forall j = 1,2,\ldots,n$$
$$e^T x = 1,$$
$$x \geq 0,$$
$$0 \leq \lambda \leq 1,$$

(FD2) min $F(\widetilde{w})$

subject to ,

$$\sum_{j=1}^{n} F(\tilde{a}_{ij}) y_j \gtrsim F(\widetilde{w}) + (1-\eta)\, F(\tilde{q}) \qquad \forall\ i = 1,2,\ldots,m$$
$$e^T y = 1,$$
$$y \geq 0,$$
$$0 \leq \eta \leq 1,$$

**Example 2.1** Consider the fuzzy matrix game with fuzzy payoff matrix

$$\begin{pmatrix} (175,180,190) & (150,156,158) \\ (80,90,100) & (175,180,190) \end{pmatrix}$$

Where ($a_l$, $a$, $a_u$) is a triangular fuzzy number $\widetilde{p_1}$ assuming that Player I and Player II have the margins $\widetilde{p_1} = \widetilde{p_2} =$(0.08, 0.10, 0.11), and $\widetilde{q_1} = \widetilde{q_2}$ = (0.14, 0.15, 0.17).

**Solution.**

To solve this game, following two crisp linear programming problems (LP) and (LD) is needed to be solve for Player I and Player II respectively:

(LP) $\quad \max V=\frac{(v_1+v_2+v_3)}{3}$

Subject to,

$$545x_1 + 270x_2 \geq 3V - (1-\lambda)*(0.29)$$
$$464x_1 + 545x_2 \geq 3V - (1-\lambda)*(0.29)$$
$$x_1 + x_2 = 1$$
$$\lambda \leq 1$$
$$x_1, x_2, \lambda \geq 0,$$

(LD) $\quad \min W=\frac{(w_1+w_2+w_3)}{3}$

Subject to,

$$545y_1 + 464y_2 \leq 3W + (1-\eta)*(0.46)$$
$$270y_1 + 545y_2 \leq 3W + (1-\eta)*(0.46)$$
$$y_1 + y_2 = 1$$
$$\eta \leq 1.$$
$$y_1, y_2, \eta \geq 0.$$

Solving the above (LP) and (LD), the values are, ($x_1^*$=0.7725, $x_2^*$= 0.2275, V = 160.91, $\lambda$∗= 0) and ($y_1^*$= 0.2275, $y_2^*$=0.7725, W = 160.65, $\eta$∗= 0) respectively. Therefore, optimal strategies for Player I and Player II are ($x_1^*$= 0.7725, $x_2^*$= 0.2275) and ($y_1^*$= 0.2275, $y_2^*$ = 0.7725) respectively.

**Remark 2.1** Here $F(\tilde{v})$ and $F(\tilde{w})$ are used which is in crisp sense instead of using $\tilde{v}$ and $\tilde{w}$ in fuzzy sense in objective function of (FP) and (FD). So it gives an average of fuzzy values instead of fuzzy value for Player I and Player II.

Therefore, another method which is based on "Pareto-optimal security strategies" used to solve Two person zero-sum matrix games with Fuzzy Pay-offs.

## 3. Pareto-optimal security strategies in matrix games with fuzzy payoffs:

If the membership function is given by a piecewise linear function then it has a finite number of pieces and therefore only a finite number of different α-cut sets is used to exactly describe the corresponding fuzzy number. Set of cuts is defined as $\gamma = \{\alpha_1, \alpha_2, \ldots\ldots, \alpha_r\} \subset [0, 1]$ with $\alpha_1 < \alpha_2 < \cdots < \alpha_{r-1} < \alpha_r = 1$ α-cut of a fuzzy number can be defined as $\tilde{a}$ by $\widetilde{a_\alpha} = [a_\alpha^1, a_\alpha^2]$ where, $a_\alpha^1$= inf ($a_\alpha$) and $a_\alpha^2$=sup ($a_\alpha$).

**Standard ranking function:**

A function f: N(R) →$R^{r*2}$ such that, f ($\tilde{a}$) = ($p_{ij}(\tilde{a})$) ∈ $R^{r*2}$ where $p_{ij}(\tilde{a}) = a_{\alpha_i}^j$ , i = 1,…, r and j = 1, 2.

**Standard fuzzy order:**

A fuzzy order, $\lesssim$ is standard if there exists a standard ranking function f such that, $\tilde{a} \lesssim \tilde{b} \Leftrightarrow$ f $(\tilde{a}) \le$ f $(\tilde{b})$ with $\tilde{a}, \tilde{b} \in N(R)$, where $\le$ is the component wise order on $R^{r*2}$.

**Security Level :**

The security level of player I for strategy $x \in S^m$ is a fuzzy number $\widetilde{\underline{v}(x)} \in N(R)$ such that

$$f(\widetilde{\underline{v}(x)}) = \inf_{y \in S^n} p_{ij}(\widetilde{E(x,y)}) \qquad i=1...r \text{ and } j=1, 2$$

The security level of player II for strategy $y \in S^n$ is a fuzzy number $\widetilde{\overline{v}(y)} \in N(R)$ such that

$$f(\widetilde{\overline{v}(y)}) = \sup_{x \in S^m} p_{ij}(\widetilde{E(x,y)}) \qquad i=1...r \text{ and } j=1, 2$$

**Pareto-optimal security strategy (POSS):**

A strategy $x* \in X$ is a Pareto-optimal security strategy (POSS) of the game $\tilde{G}$ using the standard fuzzy order for player I iff there is no $x \in X$ such that

$$\widetilde{\underline{v}(x*)} \lesssim \widetilde{\underline{v}(x*)}$$

Similarly, one can define POSS for player II.

**FLP and MLP Equivalence:**

(FLP) max $\tilde{v}$

Subject to,

$$x^T\tilde{A} \gtrsim \tilde{v}$$
$$x \in S^m,$$

Where $\tilde{c}$ and $\tilde{b}$ are fuzzy vectors and $\tilde{A}$ is a matrix with fuzzy entries.

(MLP) max $(v_{ij})$ $i = 1\ldots r$, and $j = 1, 2$

Subject to,

$$x^T(\tilde{A}^j_{\alpha_i}) \ge v_{ij} * e^T \qquad \forall\, i, j$$
$$x \in S^m$$

**Example 3.1**

Consider the fuzzy matrix game with fuzzy payoff matrix

$$\begin{pmatrix} (175,180,190) & (150,156,158) \\ (80,90,100) & (175,180,190) \end{pmatrix}$$

Where $(a_l, a, a_u)$ is a triangular fuzzy number.

**Solution.**

Assume that player I wishes to solve the above game using a standard fuzzy order with set of cuts $\gamma = \{0, 1\}$

$$A_{\alpha_i} = [(a - a_l)\,\alpha_i + a_l, -(a_u - a)\,\alpha_i + a_u]$$

Take $\alpha_1=0$ and $\alpha_2=1$ to get,

$$A^1_{\alpha_1} = \begin{pmatrix} 175 & 150 \\ 80 & 175 \end{pmatrix},\ A^2_{\alpha_1} = \begin{pmatrix} 190 & 158 \\ 100 & 190 \end{pmatrix},\ A^1_{\alpha_2} = \begin{pmatrix} 180 & 156 \\ 90 & 180 \end{pmatrix},\ A^2_{\alpha_2} = \begin{pmatrix} 180 & 156 \\ 90 & 180 \end{pmatrix}$$

Use (FLP) and (MLP) equivalence to get the following (MLP)

(MLP) Max $v_{11}, v_{12}, v_{21}$

Subject to,

$$175x_1 + 80x_2 \geq v_{11} ,$$
$$150x_1 + 175x_2 \geq v_{11} ,$$
$$190x_1 + 100x_2 \geq v_{12} ,$$
$$158x_1 + 190x_2 \geq v_{12} ,$$
$$180x_1 + 90x_2 \geq v_{21} ,$$
$$156x_1 + 180x_2 \geq v_{21} ,$$
$$x_1 + x_2 = 1,$$
$$x_1, x_2 \geq 0$$

After solving this (MLP) Poss strategy for player I is $x^*$= (0.7916, 0, 0.2084) and its security level as $\widetilde{v}^*$= (161, 6.79, 3.67).

## 4. Two person zero-sum matrix game with I-fuzzy goals:

**I-fuzzy set (Intuitionistic fuzzy set):**

Let X be a universal set. An I-fuzzy set A in X is described by A = $\{\langle x, \mu_A(x), \nu_A(x)\rangle | x \in X\}$, where $\mu_A : X \rightarrow [0, 1]$ and $\nu_A : X \rightarrow [0, 1]$ define respectively the degree of belongingness and the degree of non-belongingness of an element x ∈ X to the set A with $0 \leq \mu_A(x) + \nu_A(x) \leq 1$.
If $\mu_A(x) + \nu_A(x) = 1$, ∀ x ∈ X, then A degenerates to a standard fuzzy set in X.

**Union and intersection:**

Let A and B be two I-fuzzy sets in X, their union and intersection are I-fuzzy sets defined respectively as

$A \cup B = \{\langle x, \max\{\mu_A(x), \mu_B(x)\}, \min\{\nu_A(x), \nu_B(x)\}\rangle | x \in X\}$ and
$A \cap B = \{\langle x, \min\{\mu_A(x), \mu_B(x)\}, \max\{\nu_A(x), \nu_B(x)\}\rangle | x \in X\}$.

**Score function:**

Let A be an I-fuzzy set in X, then the score function is defined as, $S(x) = \mu_A(x) - \nu_A(x)$ ∀ x ∈ X.

### 4.1 The pessimistic approach to defining the I-fuzzy statement $x \text{ (I F)} \gtrsim a$:

Let $\{\langle x, \mu(x), \nu(x)\rangle | x \in X\}$ be an I-fuzzy set.

In this approach the decision maker takes extra cautious for acceptance. That is, even if the degree of rejection of x is zero, the decision maker is not accepting totally. To represent this situation, assume that the tolerances are p and q with 0 < q < p, then the membership function is given as

$$\mu(x) = \begin{cases} 1 & x \geq a \\ 1 - \frac{a-x}{p} & a - p < x < a \\ 0 & x \leq a - p \end{cases}$$

and the non-membership function is given as

$$\nu(x) = \begin{cases} 1 & x \leq a - p \\ 1 - \frac{x-a+p}{q} & a - p < x < a - p + q \\ 0 & x \geq a - p + q \end{cases}$$

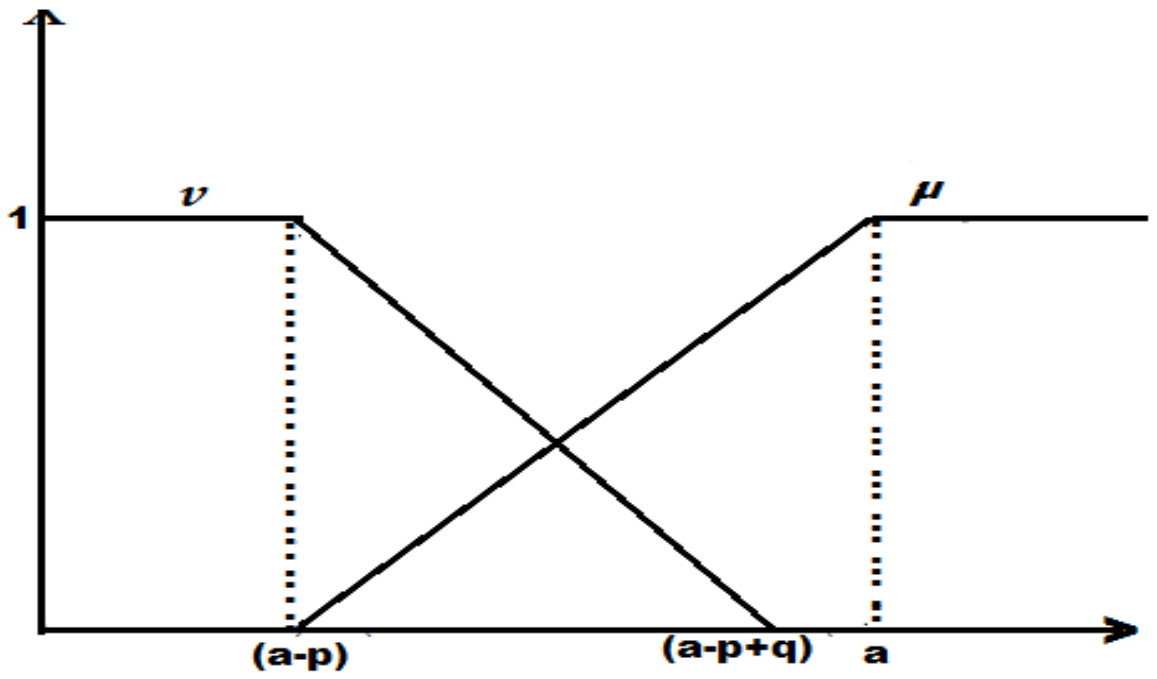

Figure 1

**Note 4.1** There is an interval [a – p + q, a] in which the non-membership degree is zero but the membership degree is not one.

**4.2 The optimistic approach to defining the I-fuzzy statement x (I F) $\gtrsim$ a:**

Let $\{\langle x, \mu(x), \nu(x)\rangle | x \in X\}$ be an I-fuzzy set. In optimistic approach the decision maker is not rejecting totally even if the degree of acceptance of x is zero. To represent this situation, there for certain tolerances p, q, the membership function is given as

$$\mu(x) = \begin{cases} 1 & x \geq a \\ 1 - \frac{a-x}{p} & a - p < x < a \\ 0 & x \leq a - p \end{cases}$$

and the non-membership function is given as

$$\nu(x) = \begin{cases} 1 & x \leq a - p - q \\ 1 - \frac{x-a+p+q}{p+q} & a - p - q < x < a \\ 0 & x \geq a \end{cases}$$

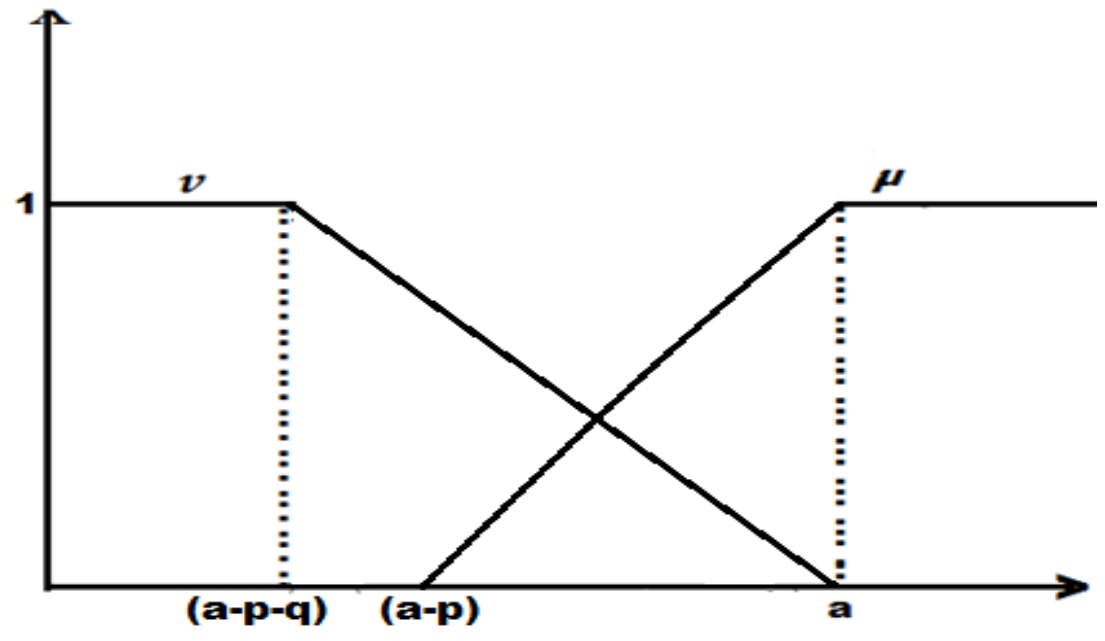

Figure 2

**Note 4.2** There is an interval [a - p - q, a - p] in which the membership degree is zero but the non-membership degree is not one.

**Remark 4.2** x (I F) $\lesssim$ a if and only if (-x) (I F) $\gtrsim$ (-a).

**Decision making in I-fuzzy environment:**

Let X be any set. Let $G_i$, for all i = 1, 2, . . . , r, be the set of r goals and $C_j$ , for all j = 1, 2, . . . ,m, be the set of m constraints, each of which can be characterized by an I-fuzzy set on X . The I-fuzzy decision is given by,

$D = (G_1 \cap G_2 \cap \cdots \cap G_r) \cap (C_1 \cap C_2 \cap \cdots \cap C_m)$ is an I-fuzzy and it is defined as ,

$$D = \{\langle x, \mu_D(x), \nu_D(x)\rangle \mid x \in X\},$$

Where, $\mu_D(x) = \min_{i,j} \{\mu_{G_i}(x), \mu_{C_j}(x)\}$ and $\nu_D(x) = \max_{i,j} \{\nu_{G_i}(x), \nu_{C_j}(x)\}$.

**I-fuzzy linear programming and duality:**

Let $R^n$ denote the n-dimensional Euclidean space. Let $c \in R^n$, $b \in R^m$, and $A \in R^{m\times n}$. Consider a general model for an I-fuzzy linear programming problem (IFP) in which the aspiration level $Z_0$ for the objective function is indicated by the decision maker.

(IFP) Find $x \in R^n$

Such that $c^T x \ (I\,F) \gtrsim Z_0$

$Ax \ (I\,F) \lesssim b$

$x \geq 0$

(IFD) Find $y \in R^m$

Such that $b^T y \ (I\,F) \lesssim W_0$

$A^T y \ (I\,F) \gtrsim c$

$y \geq 0$

Where $W_0$ is an aspiration level for the dual objective.

**Duality under pessimistic approach:**

Consider the I-fuzzy primal problem (IFP). Let $p_i$, $q_i$, $0 < q_i < p_i$, i = 0, 1, . . . , m, denote the tolerances associated respectively with the acceptance and the rejection of m + 1 constraints in (IFP). Also, let α, β denote the minimal degree of acceptance and the maximal degree of rejection respectively of the m + 1 constraints of the primal problem (IFP). The I-fuzzy optimization problem (IFP) under pessimistic scenario is equivalent to the following crisp optimization problem

(IFPC) max $\alpha - \beta$

Subject to

$$(1 - \alpha)\, p_0 + c^T x - Z_0 \geq 0$$
$$(1 - \alpha)\, p_i - A_i x + b_i \geq 0 \qquad \text{for all } i = 1, \ldots, m$$
$$(1 - \beta)\, q_0 - c^T x + (Z_0 - p_0) \leq 0$$
$$(1 - \beta)\, q_i + A_i x - (b_i + p_i) \leq 0 \qquad \text{for all } i = 1, \ldots, m$$
$$\alpha + \beta \leq 1$$
$$\alpha \geq \beta \geq 0$$
$$x \geq 0.$$

Next consider the I-fuzzy dual problem (IFD). Let $s_j$, $t_j$ , $0 < t_j < s_j$, j = 0, 1, . . . , n, be the tolerances associated respectively with the acceptance and rejection of the n + 1 constraints in (IFD). Let δ and η respectively denote the minimal degree of acceptance and the maximal degree of rejection of n + 1 constraints of the dual problem (IFD). Analogous to the discussion above, the I-fuzzy optimization problem (IFD) is equivalent to the following crisp optimization problem

(IFDC) max $\delta - \eta$

Subject to

$$(1 - \delta)\, s_0 - b^T y + W_0 \geq 0$$
$$(1 - \delta)\, s_j + A^T_j y - c_j \geq 0 \qquad \text{for all } j = 1, \ldots, n$$
$$(1 - \eta)\, t_0 + b^T y - (W_0 + s_0) \leq 0$$

$$(1-\eta)\, t_i - A^T_j\, y + (c_j - s_j) \le 0 \qquad \text{for all } j = 1, \dots, n$$
$$\delta + \eta \le 1,$$
$$\delta \ge \eta \ge 0,$$
$$w \ge 0.$$

**Remark.** The pair (IFPC) and (IFDC) is termed as I-fuzzy primal-dual pair.

### 4.3 Two person zero-sum matrix game with I-fuzzy goals:

Let $A \in R^{m\times n}$ be a real matrix. Let $U_0$ and $V_0$ be the aspiration levels for player I and II respectively, then the two person zero-sum matrix game with I-fuzzy goals is defined as,

$$IFG = (S^m, S^n, A, U_0, (I\,F)\gtrsim, p_0, q_0\,;\, V_0, (I\,F)\lesssim, s_0, t_0).$$

Where $S^m$ and $S^n$ are strategies spaces for player I and II respectively.$p_0$ and $q_0$ respectively are the tolerances pre-specified by player I for accepting and rejecting the aspiration level $U_0$. Similarly, $s_0$ and $t_0$ respectively are the tolerances pre-specified by player II for accepting and rejecting the aspiration level $V_0$.

**Solution of the I-fuzzy matrix game IFG:**

An ordered pair $(\bar{x}, \bar{y}) \in S^m \times S^n$ is said to be a solution of the two person zero-sum matrix game with fuzzy goals IFG if

$$\bar{x}^T A y\ (I\,F)\gtrsim U_0 \qquad \text{for all } y \in S^n,$$

and

$$x^T A \bar{y}\ \ (I\,F)\lesssim V_0 \qquad \text{for all } x \in S^m$$

Since $S^m$ and $S^n$ are convex polytopes, it is sufficient to consider only the extreme points of $S^m$ and $S^n$. Therefore, the problem of finding a solution of IFG is equivalent to solving the following two I-fuzzy linear programming problems

(IFG1) Find $x \in S^m$

such that $\sum_{i=1}^{m} a_{ij}x_i\ (I\,F) \gtrsim U_0, \qquad \forall\, j=1, 2, \dots, n$

$$\sum_{i=1}^{m} x_i = 1$$
$$x \ge 0$$

(IFG2) Find $y \in S^n$

such that $\sum_{j=1}^{n} a_{ij}y_j\ \ (I\,F) \lesssim V_0, \qquad \forall\, i=1, 2, \dots, m$

$$\sum_{j=1}^{n} y_j = 1$$
$$y \ge 0$$

Here (IFG1) is the I-fuzzy linear programming problem for player I and (IFG2) is the I-fuzzy linear programming problem for player II. Assuming that both players have pessimistic view points for their I-fuzzy aspiration levels. In other words, it amounts to saying that $U_0 - p_0 < U_0 - p_0 + q_0 < U_0$ and
$V_0 < V_0 + s_0 - t_0 < V_0 + s_0$, i.e., $0 < q_0 < p_0$ and $0 < t_0 < s_0$.

Let $\alpha$ denote the minimal degree of acceptance and $\beta$ denote the maximal degree of rejection of the constraints of (IFG1). Similarly let $\delta$ denote the minimal degree of acceptances and $\eta$ denote the maximal degree of rejection of the constraints of (IFG2). Then the two I-fuzzy linear programming problems are respectively equivalent to the following two crisp optimization problems.

(CFP1) Max $\alpha - \beta$

Subject to, $(1-\alpha)\, p_0 + A_j^T\, x - U_0 \ge 0 \qquad j = 1, 2, \dots, n$

$$(1-\beta)\, q_0 - A_j^T\, x + (U_0 - p_0) \le 0 \qquad j= 1, 2, \dots, n$$
$$x \ge 0, \qquad \sum_{i=1}^{m} x_i = 1$$
$$\alpha \ge \beta \ge 0,\ \ \alpha + \beta \le 1.$$

(CFP2) Max $\delta - \eta$

Subject to $(1 - \delta)\, s_0 - A_i\, y + V_0 \geq 0 \quad i = 1, . . . , m$

$(1 - \eta)\, t_0 + A_{i\cdot}\, y - (V_0 + s_0) \leq 0 \quad i = 1, . . . , m$

$y \geq 0, \quad \sum_{j=1}^{n} y_j = 1$

$\delta \geq \eta \geq 0, \; \delta + \eta \leq 1.$

Here $A_j$ and $A_i$ denote the *j* th column and the *i* th row of A, respectively.

## Appendix A.

### A.1 MATLAB code for two person zero-sum matrix game with fuzzy goals:

```
A=input ('Enter the pay-offs matrix A = ');
v0=input ('Enter the aspiration level of Player 1st  v0 = ');
w0=input ('Enter the aspiration level of Player 2nd  w0 = ');
p0=input ('Enter the tolerance error of Player 1st   p0 = ');
q0=input ('Enter the tolerance error of Player 2nd   q0 = ');
 % for player 1st
fprintf ('\n \n\n');
fprintf ('(FLP)');
x=sym ('x',[length(A(:,1)),1]);
syms v;
fprintf (' max ');
disp (v);
fprintf ('subject to;  \n');
disp(-A'*x+p0*v <=p0-v0);
disp(v<=1);
disp(sum(x)==1);
disp(x>=0);
disp(v>=0);
%for player 2nd
fprintf ('\n \n\n');
fprintf ('(FLD)');
y=sym('y',[length(A(1,:)),1]);
syms w;
fprintf(' max ');
disp(w);
```

```
fprintf('subject to;  \n');
disp(A*y+q0*w <=w0+q0);
disp(w<=1);
disp(sum(y)==1);
disp(y>=0);
disp(w>=0);
```

**A.2 MATLAB code for two person zero-sum matrix game with fuzzy pay-offs:**

```
m=input('no. of rows in matrix:m = ');
n=input('no. of columns in matrix:n = ');
disp('Enter the elements of fuzzy matrix in [al a au] form ');
for i=1:m
for j=1:n
a=input('');
A(i,j)=sum(a);
end
end
disp('Matrix is..\n');
A
disp('Enter the tollerence for player 1st in [pl p pu] form ');
p=input('p = ');
p=sum(p)
disp('Enter the tollerence for player 2nd in [ql q qu] form ');
q=input('q = ');
q=sum(q)
% for player 1st
disp('V=(vl+v+vu)\3');
fprintf('\n \n\n');
fprintf('(LP)');
x=sym('x',[m,1]);
syms V;
fprintf(' max ');
disp(V);
fprintf('subject to;  \n');
syms c;
disp(-A'*x+3*V+c*p <=p);
disp(sum(x)==1);
disp(c<=1);
disp(x>=0);
disp(c>=0);
%for player 2nd
disp('W=(wl+w+wu)\3');
fprintf('\n \n\n');
fprintf('(LD)');
y=sym('y',[n,1]);
syms W;
fprintf(' min ');
disp(W);
fprintf('subject to;  \n');
syms d;
```

```
disp(-3*W+A*y+q*d <=q );
disp(d<=1);
disp(sum(y)==1);
disp(y>=0);
disp(d>=0);
%solution of (FLP)
f1=-[1;zeros(m+1,1)];
B1=[[3*ones(n,1);1],[p*ones(n,1);1],[-A';zeros(1,m)]];
b1=[(p)*ones(n,1);1];
B1eq=[0,0,ones(1,m)];
b1eq=1;
lb1=zeros(m+2,1)';
X=linprog(f1,B1,b1,B1eq,b1eq,lb1);
fprintf('\n FOR PLAYER 1st \n')
fprintf('\n The optimal value of (FLP) is \n %f \n',X(1,1));
fprintf('\n The optimal solution of (FLP) is \n  ');
disp(X(2:m+2,1));
%solution of (FLD)
f2=-[1;zeros(n+1,1)];
B2=[[(-3)*ones(m,1);1],[q*ones(m,1);1],[A;zeros(1,n)]];
b2=[q*ones(m,1);1];
B2eq=[0,0,ones(1,n)];
b2eq=1;
lb2=zeros(n+2,1)';
Y=linprog(f2,B2,b2,B2eq,b2eq,lb2);
fprintf('\n FOR PLAYER 2st \n')
fprintf('\n The optimal value of (FLD) is \n %f \n',Y(1,1));
fprintf('\n The optimal solution of (FLD) is \n  ');
disp(Y(2:n+2,1));
```

1) **MATLAB code for two person zero-sum matrix game with fuzzy pay-offs using multi-objective linear programming:**

```
m=input('no. of rows in matrix     :m = ');
n=Input('no. of columns in matrix  :n = ');
r=input('Enter no of cuts          :r = ');
disp('Enter the elements of fuzzy matrix in(al a au) form ');
for  c=1:r
disp('Enter the value of alpha');
alpha=input('');
for i=1:m
for j=1:n
al=input('al:');
a=input('a:');
au=input('au:');
C(i,j)=(a-al)*alpha+al;
D(i,j)= (-(au-a)*alpha)+au;
end
end
C
D
```

```
end
% % for player 1st
fprintf('\n \n\n');
fprintf('(MLP)');
x=sym('x',[m,1]);
fprintf(' max v11 ,v12 ,v21');
fprintf('subject to;  \n');
for k=1:r
disp(-C'*x <=-v(k,k));
disp(-D'*x <=-v(k,k+1));
end
disp(sum(x)==1);
disp(x>=0);
%for player 2nd
disp('W=(wl+w+wu)\3');
fprintf('\n \n\n');
fprintf('(LD)');
y=sym('y',[n,1]);
syms d;
fprintf(' min ');
disp(W);
fprintf('subject to;  \n');
disp(A*y-3*W+q*d <=d );
disp(d<=1);
disp(sum(y)==1);
disp(y>=0);
disp(d>=0);
```